\documentclass[a4paper, 12pt]{article}
 \usepackage[utf8]{inputenc} 
\usepackage{makeidx}
\usepackage{enumerate, theorem}
\usepackage{amsmath, amsfonts, amssymb}
\usepackage[height=22.5cm, width=15cm]{geometry}
\usepackage[all]{xy}
\usepackage{mathrsfs, yfonts}
\usepackage{graphicx}
\DeclareMathOperator{\C}{\mathbb{C}}
\DeclareMathOperator{\codim}{codim}

\newcommand{\parag}[1]{\paragraph{\sc{#1.}} }

\newtheorem{thm}{Theorem}[subsection]

\newtheorem{cor}[thm]{Corollary}
\newtheorem{prop}[thm]{Proposition}
\newtheorem{lemma}[thm]{Lemma}

\begin{document}
\author{Daniel  Barlet\footnote{Institut Elie Cartan, G\'eom\`{e}trie,\newline
Universit\'e de Lorraine, CNRS UMR 7502   and  Institut Universitaire de France.} \ and J\'on  Magn\'usson\footnote{Department of Mathematics, School of Engineering and Physical Sciences, University of Iceland.}.}
\title{Stable properties under weakly geometrically flat maps}
\maketitle

\parag{Abstract}In this note we show that a weakly geometrically flat map $\pi\colon M\rightarrow N$ between  pure dimensional complex spaces has the local lifting property for cycles.  From this result we also deduce that, under these hypotheses, several properties of $M$ are transferred to $N$.

\section{Lifting cycles locally}

An equidimensional holomorphic map $\pi\colon M\rightarrow N$ between two pure dimensional complex spaces with $q := \dim M - \dim N$,  is called {\bf  weakly geometrically flat} (or geometrically flat in the weak sense)  if it is surjective and  there exists an analytic family $(F_y)_{y \in N}$ of $q$-cycles in $M$ such that, for any $y \in N$ we have the equality $\vert F_y \vert = \pi^{-1}(y)$.\\
If, moreover, for very general $y \in N$ the cycle $F_y$ is a reduced cycle we say that $\pi$ is {\bf geometrically flat} (or geometrically flat in the strong sense if we want to avoid confusion).\\

\parag{Remarks} \begin{enumerate}[(i)]
\item 
Weak geometric flatness is stable by base change $f\colon P\rightarrow N$, where $f$ is a holomorphic map and $P$ is a pure dimensional reduced complex space. On the contrary, this is not true, in general, for strong geometric flatness when $P$ has at least one irreducible component whose image has empty interior in $N$.
\item 
When $N$ is normal and $M, N$ are pure dimensional, every equidimensional  holomorphic map $\pi\colon M \rightarrow N$ is geometrically flat (in the strong sense).
\end{enumerate}

\begin{prop}\label{Lifting} 
Let $\pi\colon M\rightarrow N$ be a geometrically flat holomorphic map in the weak sense between reduced complex spaces of pure dimensions $m$ and $n$ and $z_0$ be a point in $M$. Then there exists an open connected neighborhood $N_0$ of $\pi(z_0)$ in $N$ and a closed analytic subset $Z$ of $\pi^{-1}(N_0)$ such that $\pi^{-1}(\pi(z_0))\cap Z = \{z_0\}$ and $\pi$ induces a proper and finite map $\pi_0\colon Z\rightarrow N_0$, which is geometrically flat in the weak sense.  
\end{prop}

\parag{Proof} Put $y_0 := \pi(z_0)$, $q := m-n$ and let $(F_y)_{y\in N}$ be the analytic family of fiber-cycles  of $\pi$. Let $E = (U,B,j)$ be a $q-$scale on $M$,  adapted to $F_{y_0}$ and centered at $z_0$.  More precisely, $j$ is a proper holomorphic embedding of an open neighborhood $M_0$ of $z_0$ in $M$ into an open neighborhood $W$
of the origin in $\C^l$ and $U$, $B$ are open relatively compact polydisks in $\C^q$, $\C^{l-q}$ such that $\bar{U}\times\bar{B}\subset W$ and 
$$
j(F_{y_0}\cap M_0)\cap(\bar U\times\partial B) = \emptyset.
$$
Moreover, $E$ can be chosen in such a way that $j(F_{y_0}\cap M_0)\cap(\{0\}\times\bar B) = \{(0,0)\}$.
Then there exists an open neighborhood $N_0$ of $y_0$ in $\pi(M_0)$ such that, for all $y$ in $N_0$, the scale $E$ is adapted to $F_y$  and ${\deg}_{E}(F_y) = {\deg}_E(F_{y_0})$.

Put $M' := M_0\cap\pi^{-1}(N_0)$. Then $j(M') = (\pi\circ j^{-1})^{-1}(N_0)$ is an analytic subset of $W$ which satisfies 
$j(M')\cap(\bar U\times\partial B) = \emptyset$.  Hence, for all $x$ in $U$, $j(M')\cap (\{x\}\times B)$ is an analytic subset of $W$ and consequently $Z_x := j^{-1}((\{x\}\times B)\cap j(M'))$ is an analytic subset of $\pi^{-1}(N_0)$. Furthermore, the induced map 
$\pi_x\colon Z_x\rightarrow N_0$ is proper, surjective and with finite fibers and we may assume that $Z_{x_0}\cap\pi^{-1}(y_0) = \{z_0\}$.
 For every $y$ in $N_0$, denote $X_y$ the multigraph that $F_y$ determines in  $U\times B$ and put $A_y := j^*(X_y\cap_{U\times B}(\{x\}\times B))$.\\ Then $(A_y)_{y\in N_0}$ is an analytic family of $0-$cycles in $Z_x$ such that $|A_y| = \pi_x^{-1}(y)$ and it follows that 
$Z := Z_0$ satisfies  the required properties. 
\hfill{$\blacksquare$}

\parag{Remarks}
\begin{enumerate}[(i)]
\item
Under the assumptions of the proposition assume that $\pi$ is geometrically flat in the strong sense. Then for very general $y$ in $N_0$, the cycle $F_y$ is reduced.\\
So the map $\pi_0\colon Z_0 \rightarrow N_0$ is geometrically flat in the strong sense (after shrinking $N_0$ if necessary).
\item
Let $\Sigma$ be an analytic subset of codimension $r$ in $M$. Then $j(\Sigma\cap M')$ is an analytic subset of codimension  at most equal to $r$ in $j(M')$  and it follows that, for very general $x$ in $U$, we have $\codim_{Z_x}(\Sigma\cap Z_x)= r$ (or empty).  So we may choose $x_0 \in U$ such that $Z_{x_0} \cap \Sigma$ has codimension at most  equal to $r$ in $Z_{x_0}$.
\hfill{$\square$}
\end{enumerate}

\begin{cor}\label{Local lift}
Let $\pi\colon M\rightarrow N$ be a geometrically flat map  in the weak sense between pure dimensional complex  spaces, $Y$ be a rational $p$-cycle in $N$ and $y_0$ be a point in $|Y|$. Then, for every $x$ in $\pi^{-1}(y_0)$, there exists an open neighborhood $N_1$ of $y_0$ in $N$ and a rational $p$-cycle $X$ in $\pi^{-1}(N_1)$ such that $x\in |X|$ and $\pi_*X = Y\cap N_1$.
\end{cor}
\parag{Proof}  We may assume that $Y$ is an irreducible analytic subset of $N$. Then the induced  map $\pi^{-1}(Y)\rightarrow Y$ is geometrically flat in the weak sense so it is enough to consider the case where $Y = N$. Let $z_0$ be a fixed point in $\pi^{-1}(y_0)$.\\
 Then, by Proposition \ref{Lifting}, there exists an open connected neighborhood $N_0$ of $y_0$ in $N$ and an analytic subset $Z$ of $\pi^{-1}(N_0)$ such that $\pi$ induces a proper and finite  map $\pi_0\colon Z \rightarrow N_0$ which is  geometrically flat in the weak sense and such that (set-theoretically) $\pi^{-1}(z_0)\cap Z = \{z_0\}$.  As the map $\pi_0$ is both  proper and open and $N_0$ is connected, $\pi_0$ is surjective and induces a surjective map from every irreducible component of $Z_0$ onto an irreducible component of $N_0$. Let $N_0'$ be the union of those irreducible components of $N_0$ which do not contain $y_0$. Thus $N_1= N_0\setminus N_0'$ is an open connected neighborhood of $y_0$ which has only finitely many irreducible components $C_1\ldots,C_l$. For every $j\in \{1,\ldots,l\}$, we choose an irreducible component   $Z_j$ of $\pi^{-1}(N_1)$ such that $\pi(Z_j) = C_j$ and denote $k_j$ the degree of the induced map $Z_j\rightarrow C_j$. Then the rational cycle 
$$
Y := \frac 1{k_1}Z_1 +\cdots + \frac 1{k_l}Z_l
$$
has the desired properties.
\hfill{$\blacksquare$}

\section{The sheaf $L^\bullet$} 

We recall here some classical equivalent definitions of $L^2$-holomorphic forms on a normal complex space $M$.\\
Remember that $M$ has {\bf weakly rational singularities}  if and only if  $\omega^n_M = L^n_M$ where $M$ has pure dimension $n$. Thanks to \cite{[K-S]} (see also \cite{[GKKP]}) this implies that $\omega^p_M  = L^p_M$ for each $p \in [0, n]$. If, moreover, $M$ is Cohen-Macaulay, then $M$ has rational singularities.\\

The following well known lemma gives two equivalent definitions of these (coherent) sheaves.

\begin{lemma}\label{L2}
Let $M$ be a normal complex space and denote $S(M)$ its singular locus. Then the following properties are equivalent for a holomorphic $p$-form $\alpha$ defined on $M \setminus  S(M)$:
\begin{enumerate}
\item 
Let $\Sigma$ is a closed analytic subset of codimension $\geq 2$ in $M$ containing the singular set $S(M)$ in $M$. Then for every  open set $M'$ in $M$ and   any proper  holomorphic embedding $f\colon X \rightarrow M'$ of a $p$-dimensional reduced complex space such that $f(X) \not\subset  \Sigma$,  the integral 
$ \int_{K\setminus f^{-1}(\Sigma)}  f^*(\alpha)\wedge \overline{f^*(\alpha)} $ is finite for any compact set $K$ in $X$.
\item 
There exists a desingularization $\tau\colon \tilde{M} \rightarrow M$ such that $\tau^*(\alpha)$ extends as a holomorphic $p$-form on $\tilde{M}$.
\end{enumerate}
\end{lemma}

\parag{Proof} First assume that property $2.$ holds true and consider a  proper  holomorphic map $f\colon X \rightarrow M'$ such that $f(X) \not\subset \Sigma$. Remark first that we may assume that  
$\Sigma = S(M)$ since on the smooth part of $M$ the sheaf $\Omega_M^p$ has the   extension  property in codimension $2$. Let $\tilde{X} \to \tilde{M'} $ be the strict transform of $\tau$ by $f$  (so 
$\tilde{X }:= X \times_{M, strict}\tilde{M}$) and let $\tilde{f}$ be the canonical projection. Then we have
$$ 
\int_{pr^{-1}(K)} \tilde{f}^*(\tau^*(\alpha))\wedge \overline{\tilde{f}^*(\tau^*(\alpha))} < +\infty 
$$
and this integral is equal to $\int_K f^*(\alpha)\wedge \overline{f^*(\alpha)} $ since the projection  $pr : \tilde{X} \to X$ is a  (proper) modification with center in  $f^{-1}(S(M))$, for any compact $K$ in $X$.\\
Conversely, if the property $1.$ holds true, consider a desingularization $\tau\colon \tilde{M} \rightarrow M$ and take a  smooth generic point $x_0$ of an irreducible component of the exceptional set $E$ in $\tilde{M}$ (the preimage of the center of the modification $\tau$). Then let $g\colon X \rightarrow \tilde{M}$ be a $p$-dimensional complex submanifold of a neighborhood of $x_0$  which is transversal to $E$\footnote{Note that if $E$ has codimension $\geq 2$   the proof is complete since $\tau^*(\alpha)$ extends holomorphically across $E$.} at the point $x_0$. Then let $f := g\circ \tau$. The property $1.$ gives the fact that $g^*(\tau^*(\alpha))$ cannot have a pole along $E$, since $x_0$ is generic, and then  $\tau^*(\alpha)$ extends holomorphically across $E$, using again the codimension $2$ extension property of the sheaf $\Omega^p_{\tilde{M}}$ since $\tilde{M}$ is a complex manifold.$\hfill \blacksquare$

\parag{Remark} 
Let $Z$ be an analytic subset of an open subset $M'$ of $M$, such that $Z \cap S(M)$ is of codimension at least $2$ in $Z$, and $j\colon Z \hookrightarrow M'$ be the canonical injection. Then, due to  characterisation $1.$ above,  we have a restriction map
$$ 
j^*(L^p_{M'}) \to L^p_Z 
$$
extending the pull-back morphism of holomorphic $p$-forms, for any $p \geq 0$.\\

\section{The Theorem}

\begin{thm}\label{facile}
Let $\pi : M \to N$ be a geometrically flat map  in the weak sense between irreducible complex spaces. Assume that $M$ has one of the following properties:
\begin{enumerate}[P.1]
\item 
The complex space  $M$ is normal.
\item   
The complex space  $M$ is normal and  every codimension 1 cycle in $M$ is $\mathbb{Q}$-Cartier.
\item   
The complex space  $M$ is normal and every holomorphic $p$-form on $M \setminus S(M)$ extends as a section of the $L_M^p$-sheaf (so $L_M^p = \omega_M^p$).
\item   
The complex space  $M$ is normal  and  the de Rham complex $0 \to \C \to \big(\omega^\bullet_M, d^\bullet\big)$ is exact in degree $p$.
\item    The complex space  $M$ is normal   and  the de Rham complex $0 \to \C \to \big(L^\bullet_M, d^\bullet)$ is locally exact in degree $p$ on $M$.

\end{enumerate}
Then the same property is also true for $N$.
\end{thm}

\parag{Proof} As all five properties are local we may assume, due to Proposition \ref{Lifting} and Remark $(ii)$ following it, that there exists an analytic subset $Z$ of dimension $n$ in $M$ satisfying the following conditions:
\begin{itemize}
\item The restriction $\pi_{\vert Z} $ is surjective, proper and geometrically flat in the weak sense.
\item $S(M) \cap Z$ is of codimension at least 2. in $Z$
\end{itemize}
Denote $k$ the degree of $\pi' := \pi_{\vert Z} $.
\smallskip

Let us begin by property $P.1$. Take a locally bounded meromorphic function $f$ in an open  neighborhood $N'$  of $y_0 \in N$. Then $\pi^*(f)$ is holomorphic on $\pi^{-1}(N')$ and induces a holomorphic function $g$ on $\pi'^{-1}(N')$. Then the trace of $g$ by  the induced map $\pi'^{-1}(N') \rightarrow N'$ is holomorphic and equal to $kf$.\\

Proof of $P.2$ for $N$. Assume that $X$ is a cycle of  codimension $1$  in $N$ and $y$ a point in $N$. Then $\pi^*(X)$ is a cycle of codimension $1$  in $M$ and we may find an open neighborhood $U$ in $M$  of $\pi^{-1}(y) \cap Z$ and a holomorphic function $f$ defining the cycle $\pi^*(X)$ in $U$. As $\pi'$ is a proper map, there exists an open neighborhood $N'$ of $y$ such that $\pi'^{-1}(N') \subset U\cap Z$. Put $Z' := \pi'^{-1}(N')$ and denote  $\pi''\colon Z' \rightarrow N'$ the restriction of $\pi'$. Then $f_{|Z'}$  defines the cycle $(\pi'')^*(X)$. Since the norm of the function $f_{\vert Z'}$ is meromorphic and locally bounded on $N'$, which is normal, it is holomorphic on $N'$ and defines the cycle $k(X \cap N')$. Hence $X$ is locally $\mathbb{Q}$-Cartier in $N$.\\

To prove property $P.3$ for $N$, consider a $p$-form $\alpha$  on $N \setminus S(N)$. Then $\pi^*(\alpha)$ extends as a global section of the sheaf $L^p_M$. Since $Z$ meets $S(M)$ in codimension at least $2$ the form $\pi'^*(\alpha)$ is holomorphic on $Z \setminus S(M)$ and is a section of $L^p_Z$ due to the characterization 1. in Lemma \ref{L2}. Then, since the map $\pi'$ is a proper map with finite fibers, using again the characterization 1 above with $\Sigma = \pi(S(M)\cap Z) \cap S(N)$ we conclude that $\alpha$ is a section of $L^p_N$ because the normality of $M$ implies the normality of $N$ (see P.1), and $S(N)$ has  codimension at least $2$  in $N$.$\hfill \blacksquare$\\

To prove property $P.4$ for $N$ we take an open neighborhood $N'$  of a point $y$ in $N$ and we consider $\alpha \in \Gamma(N', \omega^p_N)$ such that $d\alpha = 0$. Then the restriction of $\alpha$ to $N' \setminus S(N)$ is a closed holomorphic $p$-form. Then $\pi^*(\alpha)$ is also a closed holomorphic $p$-form on $ \pi^{-1}(N'\setminus S(N))$. Since $\pi^{-1}(S(N))$ and  $S(M)$ have codimension at least equal to $2$, $\pi^*(\alpha)$  extends as a $d$-closed section $\beta$ of $\omega^p_M$ on $\pi^{-1}(N')$. Then we may find an open neighborhood $U$ in $M$  of $\pi^{-1}(y)\cap Z$ (which is a finite set) and a section $\beta \in \Gamma(U, \omega^{p-1}_M)$ such that 
$d\beta = \pi^*(\alpha)$ on $U$. Then $\beta$ is a $(p-1)$-holomorphic form on $U \setminus S(M)$ and then induces a $(p-1)$-holomorphic form on $Z \setminus (Z \cap S(M))$. Since $(Z \cap S(M))$ has codimension at least $2$ in $Z$ we obtain a section on $Z$ of $\omega^{p-1}_Z$ (the sheaf $\omega_Z^\bullet$ has the extension property in codimension $\geq 2$) and we have $d\beta_{\vert Z} = \pi^*(\alpha)_{\vert Z}$ since the sheaf $\omega_Z^\bullet$ has no torsion. The direct  image by $\pi'$ gives a section $\gamma := \pi'_*(\beta_{\vert Z})$  on $N'$ which satisfies $d\gamma = k \alpha$ where $k$ is the degree of $\pi'$. This proves property $P.4$ for $N$. \\ 
The proof of $P.5$ is analogous.$\hfill \blacksquare$

 \parag{Remark} In the paper \cite{[K-S]} it is proved (see Theorem 1.4\footnote{In fact they prove more: if $ \omega^k_M = L^k_M$ for some $k \leq n$ then $\omega^p_M = L^p_M$ for any $p \in [0, k]$.}), generalizing the result of \cite{[GKKP]},   that $\omega^n_M = L^n_M$ for a normal complex space implies $\omega^p_M = L^p_M$ for any $p \in [0, n]$.\\
 An interesting characterization of $L^p$ is also given in \cite{[K-S]} (see Theorem 1.1): \\
 A section of $\omega^p$ is in $L^p$ is for any $\beta \in \Omega_M^{n-p}$ and any $\gamma \in \Omega^{n-p-1}_M$ the forms $\alpha \wedge \beta$ and $d\alpha\wedge \gamma$ are in $L^n$. $\hfill \square$\\

\end{document}